\documentclass[11pt,reno]{amsart}  
\usepackage{amssymb}
\usepackage{mathtools}
\usepackage{enumitem}
\usepackage{mathrsfs}
\usepackage[english]{babel}
\usepackage[all]{xy}
\usepackage{hyperref}
\usepackage{marginnote}
\usepackage{stmaryrd}
\usepackage{fullpage}

\RequirePackage{mathrsfs}

\newtheorem{theorem}{Theorem}[section]
\newtheorem{lemma}[theorem]{Lemma}
\newtheorem{proposition}[theorem]{Proposition}
\newtheorem{corollary}[theorem]{Corollary}
\newtheorem{conjecture}[theorem]{Conjecture}

\theoremstyle{definition}
\newtheorem*{ack}{Acknowledgements}
\newtheorem*{con}{Conventions}

\newtheorem{example}[theorem]{Example}
\newtheorem{definition}[theorem]{Definition}

\numberwithin{equation}{section} \numberwithin{figure}{section}

\DeclareMathOperator{\Spec}{Spec}

\DeclareMathOperator{\Rad}{Rad}
\DeclareMathOperator{\der}{der}
\DeclareMathOperator{\Res}{Res}

\newcommand{\GG}{\mathbb{G}}
\newcommand{\Gm}{\GG_{\mathrm{m}}}
\newcommand{\Gmn}[1]{\GG_{\mathrm{m},{#1}}}
\newcommand{\adeles}{\mathbb{A}}

\title{Hilbert irreducibility for integral points on punctured linear algebraic groups}
\author{Cedric Luger}
\address{Cedric Luger\\
Institut f\"{u}r Mathematik\\
Johannes Gutenberg-Universit\"{a}t Mainz\\
Staudingerweg 9, 55099 Mainz\\
Germany.}
\email{celuger@uni-mainz.de}



\begin{document}

\begin{abstract}
    Let $K$ be a number field, let $X$ be a smooth integral variety over $K$,
    and assume that there exists a finite set of finite places
    $S$ of $K$ such that the $S$-integral points on $X$ are dense.
    Then the combined conjectures of Campana and Corvaja--Zannier predict that, for every closed subscheme $Z$ of $X$ of codimension at least two, there exists a finite extension $L$ of $K$ and a finite set of finite places $T$ of $L$ such that the $T$-integral points on $(X\setminus Z)_L$ are not strongly thin.
    The main goal of the present paper is to show that this property holds for all connected linear algebraic groups.
    Our result builds mainly on recent work on a Hilbert irreducibility type theorem for connected algebraic groups, the purity of strong approximation for semi-simple simply connected quasi-split linear algebraic groups, and the relation between integral strong approximation and the Hilbert property.
\end{abstract}

\maketitle
\thispagestyle{empty}

\section{Introduction}

Integral points on algebraic varieties have been the object of extensive study in Diophantine geometry. A central theme is to understand how the geometry of a variety constrains the existence and distribution of its integral points.

Beyond the question of existence, one may ask how ``large'' a set of integral points can be from an arithmetic perspective.
In this direction, we follow \cite{CZHP} (see also \cite[Definition~1.2]{CDJLZ} and \cite[Definition~1.3]{LugerProducts}) and study the weak Hilbert property, which is a refined version of the usual Hilbert property, as defined in \cite[\S 3]{SerreTopicsGalois}.
To give its definition, we call
a morphism $Y\to X$ of integral varieties a \emph{finite cover} if it is finite surjective separable.
If $X$ and $Y$ are normal, a finite cover $Y \to X$ is a \emph{ramified cover} if it is not étale.

\begin{definition} 
Let $X$ be a normal integral variety over a field  $k$.
A subset $\Omega\subseteq X(k)$  is \emph{strongly thin} in $X$ if there is an integer $n\geq 0$ and a finite collection $(\pi_i \colon Y_i\to X)_{i=1}^n$ of ramified covers such that  
\[
    \Omega \setminus \bigcup_{i=1}^n \pi_i (Y_i(k))
\]
is not Zariski dense in $X$.
\end{definition}

\begin{definition}
Let $S$ be a finite set of finite places of a number field $K$ and
let $X$ be a smooth integral variety over $K$.
Let $\mathcal{X}$ be an integral scheme that is dominant, separated, and of finite type over $\mathcal{O}_{K,S}$
with generic fiber $\mathcal{X}_K \cong X$.
We say that $\mathcal{X}$ has the \emph{weak Hilbert property over $\mathcal{O}_{K,S}$} if $\mathcal{X}(\mathcal{O}_{K,S})$ is not a strongly thin subset of $X=\mathcal{X}_K$.
We say that $X$ has the \emph{weak Hilbert property over $\mathcal{O}_{K,S}$} if there exists
an integral scheme $\tilde{\mathcal{X}}$ that is dominant, separated, and of finite type over $\mathcal{O}_{K,S}$,
has generic fiber $\tilde{\mathcal{X}}_K \cong X$, and has the weak Hilbert property over $\mathcal{O}_{K,S}$.
\end{definition}

For example, any smooth projective integral variety over a number field $K$ with the Hilbert property has the weak Hilbert property over $\mathcal{O}_{K,S}$ for any finite set of finite places $S$ of $K$.
This includes rational varieties by Hilbert's irreducibility theorem, but also certain non-rational del Pezzo surfaces \cite{Demeio1, Demeio2}, certain  K3 surfaces \cite{CZHP}, and abelian varieties with a dense set of rational points~\cite{CDJLZ}.

For $S$ a scheme, a \emph{torus over $S$} is a smooth affine group scheme $T$ over $S$ with geometrically connected fibers such that $T$ is fppf locally isomorphic to $\mathbb{G}_m^n$ for some integer $n\geq 0$.
In \cite{ZannierDuke} Zannier proved the weak Hilbert property for tori:
given a suitable ramified covering, he constructed a torsion point which does not lift along the ramified covering to a rational point and then used Kummer theory to construct a dense subset of integral points which do not lift.

Corvaja--Zannier asked whether the density of integral points on a variety implies the weak Hilbert property, possibly up to an extension of the base field:

\begin{conjecture}[{\cite[Question-Conjecture 1bis]{CZHP}}]\label{Conj:CZ_1bis}
Let $S$ be a finite set of finite places of a number field $K$ and
let $X$ be a smooth integral variety over $K$.
Let $\mathcal{X}$ be an integral scheme that is dominant separated of finite type over $\mathcal{O}_{K,S}$
such that $\mathcal{X}_K \cong X$.
If $\mathcal{X}(\mathcal{O}_{K,S})$ is Zariski dense in $X$, then there exists a finite field extension $L/K$ and a finite set of finite places $T$ of $L$
such that $X_L$ satisfies the weak Hilbert property over~$\mathcal{O}_{L,T}$.
\end{conjecture}

This question has been answered positively in numerous cases. For examples, see~\cite{CocciaCubic,   Demeio1, Demeio2,  GvirtzChenMezzedimi,    LugerProducts}.
The main result of this paper is a positive answer to Conjecture~\ref{Conj:CZ_1bis} for ``punctured'' linear algebraic groups, i.e., the complement of a closed subscheme of codimension at least~$2$ in a linear algebraic group:

\begin{theorem}\label{Thm:WHP_Punctured_LinAlgGrps}
    Let $G$ be a connected linear algebraic group over $K$ and $Z \subseteq G$ a closed subscheme of codimension at least~$2$.
    Then there exists a finite field extension~$L/K$ and a finite set of finite places $T$ of $L$ such that $(G\setminus Z)_L$ satisfies the weak Hilbert property over $\mathcal{O}_{L,T}$.
\end{theorem}

Theorem~\ref{Thm:WHP_Punctured_LinAlgGrps} was first obtained and explained to us by Ariyan Javanpeykar and Julian Demeio in the case of tori. The proof for the toric case has three ingredients:
 
\begin{itemize}
    \item A product theorem for arithmetic schemes with the weak Hilbert property, established in our previous paper \cite{LugerProducts} (see also Theorem~\ref{Thm:ProductsNotStronglyThin}),
    \item Zannier's \emph{strong} version of Hilbert irreducibility for tori \cite{ZannierDuke}, and
    \item the approximation argument to prove  potential density of integral points on punctured tori by Hassett--Tschinkel \cite{HT01}, which they attribute to McKinnon.  
\end{itemize}
We expand the arguments of Demeio and Javanpeykar to all linear algebraic groups by using a recently established strong version of Hilbert irreducibility for connected algebraic groups \cite{LugerProducts}, which builds on the work of Liu \cite{FeiLiu} and Corvaja--Demeio--Javanpeykar--Lombardo--Zannier \cite{CDJLZ}.
Moreover, we use a combination of the work of Cao--Xu \cite{CaoXu}, Cao--Liang--Xu \cite{CaoLiangXu}, and Wei \cite{Wei}
on strong approximation
with the work of Nakahara--Streeter \cite{NakaharaStreeter} and Moerman \cite{Moerman}
on integral strong approximation and the Hilbert property.
To properly combine these results, we introduce the property ``(OPS)" (Definition~\ref{Def:OPS}) for schemes.
By definition, a scheme with this property that has an integral point inside the complement of closed subscheme of a given codimension admits a non-strongly-thin subset of integral points in this complement.

This paper is motivated by a conjecture due to Campana \cite[Conjecture~9.20]{CampanaOrbifolds} which predicts that a variety over a number field admits a dense set of integral points if and only if it is ``special'' in the sense of \cite[Definition~2.1]{CampanaOrbifolds} (see also \cite[Definition~1.2]{BartschJavanpeykarLevin}),
and by the fact that smooth special varieties remain special after removing a closed subset of codimension at least~$2$ (see \cite[Theorem~1.21]{BartschJavanpeykarLevin}).
In particular, since Campana's class of smooth special varieties is invariant under removing closed subsets of codimension at least~$2$, the same should hold for any property conjecturally equivalent to being special, such as ``having a dense entire curve'' or ``having a potentially dense set of integral points''.
Thus, as Corvaja--Zannier conjectured that the weak Hilbert property holds for any smooth integral variety with a dense set of integral points (Conjecture~\ref{Conj:CZ_1bis}), the conjectures of Campana and Corvaja--Zannier predict the following ``puncturing conjecture'':

\begin{conjecture}\label{conj}
    Let $X$ be a smooth integral variety over a number field $K$. Then the following are equivalent.
    \begin{enumerate}[label=(\arabic*)]
    \item\label{Item:Conj_X_special} The variety $X$ is special.
    \item\label{Item:Conj_X_denseIntegralPoints} There is a number field $L/K$, a finite set of finite places $T$ of $L$, and a  model $\mathcal{X}$ for $X_L$ over $\mathcal{O}_{L,T}$ such that $\mathcal{X}(\mathcal{O}_{L,T})$ is dense.
    \item\label{Item:Conj_XminusZ_denseIntegralPoints} For every closed subscheme $Z\subseteq X$ of codimension at least two,  there is a number field $L/K$, a finite set of finite places $T$ of $L$, and a  model $\mathcal{U}$ for $(X\setminus Z)_L$ over~$\mathcal{O}_{L,T}$ such that $\mathcal{U}(\mathcal{O}_{L,T})$ is dense.
    \item\label{Item:Conj_X_IWHP} There is a number field $L$ and a finite set of finite places $T$ of $L$ such that $X_L$ satisfies the  weak Hilbert property over $\mathcal{O}_{L,T}$.
    \item\label{Item:Conj_XminusZ_IWHP} For every closed subscheme $Z\subseteq X$ of codimension at least two,  there is a number field $L$ and a finite set of finite places $T$ of $L$ such that $(X\setminus Z)_L$ satisfies the  weak Hilbert property over $\mathcal{O}_{L,T}$.  
    \end{enumerate}
\end{conjecture}
Note that \ref{Item:Conj_XminusZ_IWHP} implies \ref{Item:Conj_X_denseIntegralPoints}-\ref{Item:Conj_X_IWHP}.
The equivalence of \ref{Item:Conj_X_special} and \ref{Item:Conj_X_denseIntegralPoints} is Campana's conjecture \cite[Conjecture~9.20]{CampanaOrbifolds}, and the equivalence of \ref{Item:Conj_X_denseIntegralPoints} and \ref{Item:Conj_X_IWHP} is Conjecture~\ref{Conj:CZ_1bis} for the variety~$X$.
The equivalence of \ref{Item:Conj_X_special}, \ref{Item:Conj_X_denseIntegralPoints} and \ref{Item:Conj_XminusZ_denseIntegralPoints} is again a special case of Campana's conjectures by the aforementioned result that $X\setminus Z$ is special.
We note that part of this equivalence was also conjectured by Hassett--Tschinkel \cite[Problem~2.13]{HT01}, albeit with the weaker assumption that $X$ has only canonical singularities, in which case there are counterexamples given in \cite[Theorem~E]{BartschJavanpeykarLevin}.
Finally, the equivalence of \ref{Item:Conj_XminusZ_denseIntegralPoints} and \ref{Item:Conj_XminusZ_IWHP} is the natural extension of Conjecture~\ref{Conj:CZ_1bis} for the variety $X\setminus Z$.  

If $G$ is a connected algebraic group, then $G$ is special, satisfies \ref{Item:Conj_X_IWHP} by \cite{LugerProducts}, and thus also satisfies \ref{Item:Conj_X_denseIntegralPoints}.
Our main result, Theorem~\ref{Thm:WHP_Punctured_LinAlgGrps}, is
that connected linear algebraic groups satisfy \ref{Item:Conj_XminusZ_IWHP}
and therefore shows that Conjecture~\ref{conj} holds for all linear algebraic groups.
However, a complete proof of Conjecture~\ref{conj} remains elusive in the case of a general connected algebraic group, e.g., an abelian variety.

\begin{con}
    Throughout this paper, we let $K$ denote a number field, $S$ a finite set of finite places of $K$, and $R=\mathcal{O}_{K,S}$ the ring of $S$-integral numbers of $K$.
    By an \emph{arithmetic scheme} over $R$ we mean a scheme that is dominant separated of finite type over $R$ with an integral generic fiber. (Note that, contrary to \cite{LugerProducts}, we do not require an arithmetic scheme to be normal.)
    If $X$ is an integral variety over $K$, by a \emph{model} of $X$ over $R$ we mean an integral arithmetic scheme $\mathcal{X}$ over $R$ with generic fiber $\mathcal{X}_K \cong X$.
\end{con}

\begin{ack}
    We gratefully acknowledge Ariyan Javanpeykar for constant support and many helpful discussions. We also thank Finn Bartsch, Manuel Blickle, and Boaz Moerman for helpful discussions.
    We are grateful to Olivier Wittenberg for explaining to us the relation between strong approximation and the integral Hilbert property, and also the relevance of the work of Cao--Liang--Xu \cite{CaoLiangXu} on strong approximation for punctured semi-simple simply connected groups.
    We are grateful to Ariyan Javanpeykar and Julian Demeio for allowing us to include and build on their initial argument for punctured tori.
    We thank the referee for carefully reading the manuscript and many valuable suggestions.
    This research was funded by the Studienstiftung des Deutschen Volkes (German Academic Scholarship Foundation). We gratefully acknowledge support by the Deutsche Forschungsgemeinschaft (DFG, German Research Foundation) – Project-ID 444845124 – TRR 326. 
\end{ack}

\section{One point suffices}

In order to to verify part \ref{Item:Conj_XminusZ_IWHP} of Conjecture~\ref{conj} for certain varieties,
we investigate arithmetic schemes for which the existence of one integral point in the complement of a closed subscheme implies the non-strongly-thinness of such points.
We will say that such arithmetic schemes satisfy the property ``(OPS)'', for ``one point suffices''. If this property is obtained after an enlargement of $S$ and preserved under further enlargements, we refer to it as ``(OPSE)'', for ``one point suffices after extension''.

\begin{definition}\label{Def:OPS}
    Let $\mathcal{X}$ be an arithmetic scheme over $R$ with $\mathcal{X}_K$ normal and
    let~$c \geq 1$ be an integer.
    We say that $\mathcal{X}$ satisfies \emph{(OPS) for codimension at least~$c$ over $R$} if
    $\mathcal{X}(R)\neq \emptyset$ and,
    for every closed subscheme $\mathcal{Z} \subseteq \mathcal{X}$ such that $\mathcal{Z}_K \subseteq \mathcal{X}_K$ is of codimension at least~$c$ and
    $(\mathcal{X}\setminus \mathcal{Z})(R) \neq \emptyset$,
    the set $(\mathcal{X}\setminus \mathcal{Z})(R)$ is not strongly thin in~$\mathcal{X}_K$.
    We say that $\mathcal{X}$ satisfies \emph{(OPSE) for codimension at least~$c$}
    if there exists a finite set of finite places $S'$ of $K$ containing $S$ such that,
    for every finite set of finite places~$S''$ of $K$ containing $S'$, the scheme
    $\mathcal{X}_{\mathcal{O}_{K,S''}}$
    satisfies (OPS) for codimension at least~$c$.
\end{definition}

In the definition of (OPS) we require $\mathcal{X}(R)\neq \emptyset$
to avoid tautologies for schemes without integral points.
Note that, if $(\mathcal{X}\setminus \mathcal{Z})(R)$ is not strongly thin in $\mathcal{X}_K$, then it is also not strongly thin in $\mathcal{X}_K \setminus \mathcal{Z}_K$, as can be seen by normalizing a given ramified cover of $\mathcal{X}_K \setminus \mathcal{Z}_K$ in $\mathcal{X}_K$, which produces a ramified cover of $\mathcal{X}_K$ by purity of the branch locus.
Even though the property (OPS) will be stable under enlargements of $S$ for all schemes that we consider in this paper, this might not be the case in general,
which is why in the definition of (OPSE) we require (OPS) to hold for all enlargements $S''$ of $S'$.
An arithmetic scheme satisfying (OPS) or (OPSE) for any codimension has a dense set of rational points, since the empty subscheme is of codimension~$\infty$.

Given a closed subscheme $\mathcal{Z} \subseteq \mathcal{X}$ and a finite set of finite places $S'$ of $K$ containing~$S$, we have
\[
    (\mathcal{X}_{\mathcal{O}_{K,S'}}\setminus \mathcal{Z}_{\mathcal{O}_{K,S'}})(\mathcal{O}_{K,S'})
    = (\mathcal{X}\setminus \mathcal{Z})_{\mathcal{O}_{K,S'}}(\mathcal{O}_{K,S'})
    = (\mathcal{X}\setminus \mathcal{Z})(\mathcal{O}_{K,S'})
\]
as subsets of $\mathcal{X}_K(K)$; we make this identification tacitly when using the property~(OPS) for $\mathcal{X}_{\mathcal{O}_{K,S'}}$.

In this section, we investigate integral points on products whose factors satisfy the property (OPSE).
We use the following special case of a product theorem given in \cite{LugerProducts}:

\begin{theorem}[{\cite[Theorem~4.4]{LugerProducts}}]\label{Thm:ProductsNotStronglyThin}
	Let $\mathcal{X}, \mathcal{Y}$ be normal quasi-projective arithmetic schemes over $R$ and define
    $X=\mathcal{X}_K, Y = \mathcal{Y}_K$.
    Let $p \colon X \times Y \to X$ denote the projection morphism.
    Let $\Sigma \subseteq (\mathcal{X} \times \mathcal{Y})(R)$ be a subset.
    If $p(\Sigma)$ is not strongly thin in $X$
    and, for every $x\in p(\Sigma)$, the set
    $\Sigma\cap (\{x\}\times Y)$ is not strongly thin in $\{x\}\times Y$,
	then $\Sigma$ is not strongly thin in $X \times Y$.
\end{theorem}

The following result shows that a punctured product of schemes satisfying (OPSE)
admits a non-strongly-thin set of integral points, up to enlarging the base ring.

\begin{theorem}\label{Thm:Products_OPSE}
    Let $N \geq 1$ be an integer, let
    $\mathcal{X} = \mathcal{X}_1 \times_R \ldots \times_R \mathcal{X}_N$
    with each $\mathcal{X}_i$ a normal quasi-projective arithmetic scheme over $R$ satisfying (OPSE) for codimension at least~$c_i$, and let $c = \max_{i=1,\ldots,N}c_i$.
    Let $\mathcal{Z} \subseteq \mathcal{X}$ be a closed subscheme such that~$\mathcal{Z}_K$ is of codimension at least~$c$ in $\mathcal{X}_K$.
    Then there exists a finite set of finite places $S'$ of~$K$ containing $S$ such that $(\mathcal{X}\setminus \mathcal{Z})(\mathcal{O}_{K,S'})$ is not strongly thin in $\mathcal{X}_K$.
\end{theorem}
\begin{proof}
Enlarging $S$ if necessary, we may assume that each $\mathcal{X}_{i,\mathcal{O}_{K,S'}}$ satisfies (OPS) for codimension at least~$c_i$ for every finite set of finite places $S'$ containing $S$.
We argue by induction on $N$.
Assume that $N=1$. By the property (OPS), $\mathcal{X}(K)$ is dense, so $\mathcal{X}\setminus \mathcal{Z}$ admits a $K$-point $p$. Let $S'\supseteq S$ be a finite set of finite places such that $p\in (\mathcal{X}\setminus \mathcal{Z})(\mathcal{O}_{K,S'})$.
Since $\mathcal{X}_{\mathcal{O}_{K,S'}}$ satisfies (OPS) for codimension at least~$c$
and~$(\mathcal{X}\setminus \mathcal{Z})(\mathcal{O}_{K,S'})$ is non-empty, it follows that $(\mathcal{X}\setminus \mathcal{Z})(\mathcal{O}_{K,S'})$ is not strongly thin in $\mathcal{X}_K$, as required.

Now assume that $N>1$ and define
$\mathcal{X}' = \mathcal{X}_1 \times_R \ldots \times_R \mathcal{X}_{N-1}$.
Let $\pi' \colon \mathcal{X} \to \mathcal{X}'$ be the projection onto $\mathcal{X}'$ and let $\pi_N \colon \mathcal{X}\to \mathcal{X}_N$ be the projection onto the last factor.
Let $\pi'_K$ and $\pi_{N,K}$ denote the base change of $\pi'$ and $\pi_{N}$ along $\Spec K \to \Spec R$.

Since the $\mathcal{X}_{i,K}$ are integral and have a dense set of $K$-rational points by the property (OPS), they are geometrically integral.
By \cite[Proposition~3.5]{CaoLiangXu} and the faithful flatness of $\pi'_K$ and $\pi_{N,K}$, the codimension condition spreads out;
more precisely,
there exist dense opens $U' \subseteq \mathcal{X}'_K$ and $U_N \subseteq \mathcal{X}_{N,K}$ such that,
for every $x' \in U'(K)$ and for every $x_N \in U_N(K)$,
the codimension of $(\pi'_K)^{-1}(x')\cap \mathcal{Z}_K$ in $(\pi'_K)^{-1}(x')$ and
the codimension of $(\pi_{N,K})^{-1}(x_N)\cap \mathcal{Z}_K$ in $(\pi_{N,K})^{-1}(x_N)$
are at least~$c$.

Since $\mathcal{X}(K)$ is dense in $\mathcal{X}$,
there exists a $K$-point
$p=(p',p_N) \in (\mathcal{X} \setminus \mathcal{Z})(K)$ such that~$p' \in U'(K)$ and $p_N \in U_N(K)$.
Let $S' \supseteq S$ be a finite set of finite places of~$K$ such that $p\in (\mathcal{X}\setminus \mathcal{\mathcal{Z}})(\mathcal{O}_{K,S'})$ and define $R' = \mathcal{O}_{K,S'}$.
Let $\pi_{N,R'}$ and $\pi'_{R'}$ be the pullback of $\pi_{N}$ and $\pi'$ along
$\Spec R' \to \Spec R$, respectively,
identify $\pi_{N,R'}^{-1}(p_N)$ with~$\mathcal{X}'_{R'}$, and write $\mathcal{Z}'\coloneqq \mathcal{Z}_{R'}\cap \pi_{N,R'}^{-1}(p_N)$.
By the choice of $p_N$, the codimension of~$\mathcal{Z}'_K$ in $\mathcal{X}'_K$ is
at least~$c$, so by the induction hypothesis,
there exists a finite set of finite places $S''$ of $K$ containing $S'$ such that
$(\mathcal{X}'\setminus \mathcal{Z}')(\mathcal{O}_{K,S''})$ is not strongly thin in~$\mathcal{X}'_K$.
Define~$R'' = \mathcal{O}_{K,S''}$ and
\[
    \Gamma \coloneqq (\mathcal{X}'_{R'} \setminus \mathcal{Z}')(R'')  \cap U'(K)
,\]
which is not strongly thin in $\mathcal{X}'_K$ since $U'\subseteq \mathcal{X}'_K$ is a dense open.
Let $s\in \Gamma$ and
let~$F_s\coloneqq (\pi'_{R''})^{-1}(s) \cong \mathcal{X}_{N,R''}$.   Let $T_s = \mathcal{Z}_{R''}\cap F_s$ be  the intersection of $F_s$ and $\mathcal{Z}_{R''}$  in $\mathcal{X}_{R''}$.

Note that $F_{s,K}\cap \mathcal{Z}_K$ has codimension at least~$c$ in $F_{s,K}$
by the choice of $U'$, and that $F_s$ satisfies (OPS) for codimension at least~$c_N$ (and thus for codimension at least~$c$).
Thus, since $(s,p_N) \in (F_s \setminus T_s)(R'')$,
we obtain that $(F_s \setminus T_s)(R'')$ is not strongly thin in $F_{s,K}$.
Define
\[
    \Sigma = \bigcup_{s\in \Gamma} (F_s \setminus T_s)(R'') \subseteq (\mathcal{X} \setminus \mathcal{Z})(R'')
.\]
Since $\Gamma = \pi'(\Sigma)$ is not strongly thin in $\mathcal{X}'_K$
and, for every $s\in \Gamma$, the set
\[\Sigma \cap (\{s\}\times \mathcal{X}_{N,K}) = (F_s\setminus T_s)(R'')\] is not strongly thin
in $F_{s,K} = \{s\} \times \mathcal{X}_{N,K}$,
it follows from Theorem~\ref{Thm:ProductsNotStronglyThin} that $\Sigma$,
and thus $(\mathcal{X}\setminus \mathcal{Z})(R'')$, is not strongly thin in $\mathcal{X}_K$, as required.
\end{proof}

One might ask whether Theorem~\ref{Thm:Products_OPSE} holds with $R'' = R$
under the assumption that~$(\mathcal{X}\setminus \mathcal{Z})(R)$ is dense in $\mathcal{X}$.
Indeed, repeating the proof under this assumption still allows to find
a point $p=(p',p_N)\in (\mathcal{X}\setminus \mathcal{Z})(R)$ satisfying both $p'\in U'(K)$ and~$p_N \in U_N(K)$.
While the fiber $\mathcal{X}'\setminus \mathcal{Z}'$ does admit an $R$-point,
it is not clear why its set of $R$-points should be dense without enlarging $R$, so one would not be able to conclude that they are not strongly thin. 
However, the implication does hold if one assumes that
$\mathcal{X}' = \mathcal{X}_1 \times_R \ldots \times_R X_{N-1}$ satisfies (OPS) for codimension at least~$c$:
As~$\mathcal{Z}'$ is of codimension at least $c$ in $\mathcal{X}$', the existence of an $R$-point would then ensure that $R$-points on $\mathcal{X}'\setminus \mathcal{Z}'$ are not strongly thin, so we may choose $R''=R$ in the proof.

Since we are not able to show that a product of schemes with (OPS) again satisfies~(OPS),
the same arguments in the proof of Theorem~\ref{Thm:Products_OPSE} only show
the following statement, in which we consider products of \emph{two} schemes:

\begin{theorem}\label{Thm:Products_OPS}
    Let $c \geq 1$ be an integer and
    let $\mathcal{X}_1$ and $\mathcal{X}_2$ be normal quasi-projective arithmetic schemes over $R$ satisfying (OPS) for codimension at least~$c$.
    Moreover, define $\mathcal{X} = \mathcal{X}_1 \times_R \mathcal{X}_2$
    and let $\mathcal{Z} \subseteq \mathcal{X}$ be a closed subscheme such that
    $\mathcal{Z}_K \subseteq \mathcal{X}_K$ is of codimension at least~$c$.
    If $(\mathcal{X}\setminus \mathcal{Z})(R)$ is dense, then it is not strongly thin in $\mathcal{X}_K$.
\end{theorem}
\begin{proof} 
We adapt the arguments of the proof of Theorem~\ref{Thm:Products_OPSE}.
Let $\pi_1$ and $\pi_2$ be the projection from $\mathcal{X}$ to $\mathcal{X}_1$ and $\mathcal{X}_2$, respectively,
and let $\pi_{1,K}$, resp. $\pi_{2,K}$, denote their base changes
along $\Spec K \to \Spec R$.

Let $i\in\{1,2\}$ and note that $\mathcal{X}_{i,K}$ is integral with a dense set of $K$-rational points by the property (OPS), so it is geometrically integral.
This shows that $\mathcal{X}_K$ is also geometrically integral.
Therefore, by \cite[Proposition~3.5]{CaoLiangXu} and the faithful flatness of $\pi_{i,K}$,
there exists a dense open $U_i \subseteq \mathcal{X}_{i,K}$ such that,
for every $x \in U_i(K)$,
the codimension of $(\pi_{i,K})^{-1}(x)\cap \mathcal{Z}_K$ in $(\pi_{i,K})^{-1}(x)$ is at least~$c$.

By the density of $(\mathcal{X}\setminus \mathcal{Z})(R)$,
there exists a point
$p=(p_1,p_2) \in (\mathcal{X} \setminus \mathcal{Z})(R)$ with~$p_1 \in U_1(K)$ and $p_2 \in U_2(K)$.
Identify $\pi_2^{-1}(p_2)$ with $\mathcal{X}_1$ and let $\mathcal{Z}_1\coloneqq \mathcal{Z}\cap \pi_2^{-1}(p_2)$.
By the choice of $p_2$, the codimension of $\mathcal{Z}_{1,K}$ in $\mathcal{X}_{1,K}$ is
at least~$c$.
Thus, since $\mathcal{X}_1$ satisfies (OPS) for codimension at least~$c$ and
$(\mathcal{X}_1\setminus \mathcal{Z}_1)(R)\ni p_1$ is non-empty, we obtain that
$(\mathcal{X}_1\setminus \mathcal{Z}_1)(R)$ is not strongly thin in $\mathcal{X}_{1,K}$.

Define
\[
    \Gamma \coloneqq (\mathcal{X}_1 \setminus \mathcal{Z}_1)(R)  \cap U_1(K)
,\]
which is not strongly thin in $\mathcal{X}_{1,K}$ since $U_1\subseteq \mathcal{X}_{1,K}$ is a dense open.
Let $s\in \Gamma$ and
let~$F_s \coloneqq (\pi_1)^{-1}(s) \cong \mathcal{X}_2$.   Let $T_s = \mathcal{Z}\cap F_s$ be  the intersection of $F_s$ and $\mathcal{Z}$  in $\mathcal{X}$.
Note that $F_{s,K}\cap \mathcal{Z}_K$ has codimension at least~$c$ in $F_{s,K}$
by the choice of $U_1$, and that $F_s$ satisfies (OPS) for codimension at least~$c$.
Thus, since $(s,p_2) \in (F_s \setminus T_s)(R)$,
we obtain that $(F_s \setminus T_s)(R)$ is not strongly thin in $F_{s,K}$.
Define
\[
    \Sigma = \bigcup_{s\in \Gamma} (F_s \setminus T_s)(R) \subseteq (\mathcal{X} \setminus \mathcal{Z})(R)
.\]
Since $\Gamma = \pi_1(\Sigma)$ is not strongly thin in $\mathcal{X}'_K$
and, for every $s\in \Gamma$, the set
\[\Sigma \cap (\{s\}\times \mathcal{X}_{2,K}) = (F_s\setminus T_s)(R)\]
is not strongly thin
in $F_{s,K} = \{s\} \times \mathcal{X}_{2,K}$,
it follows from Theorem~\ref{Thm:ProductsNotStronglyThin} that $\Sigma$,
and thus $(\mathcal{X}\setminus \mathcal{Z})(R)$, is not strongly thin in $\mathcal{X}_K$, as required.
\end{proof}

\section{The case of connected commutative group schemes}

Recall that $R = \mathcal{O}_{K,S}$.
In this section we prove that connected commutative group schemes over $R$ satisfy (OPS) under suitable conditions.
Its main result is the following proposition.

\begin{proposition}\label{Prop:HT_general}
    Let $\mathcal{G}$ be a connected commutative finite type group scheme over~$R$ and let
    \[
        0 \longrightarrow H \longrightarrow \mathcal{G}_K \longrightarrow A \longrightarrow 0
    \]
    be an exact sequence of connected algebraic groups over $K$ with $H$ linear and $A$ an abelian variety.
    Assume that $\mathcal{G}(R)$ is dense in $\mathcal{G}_K$
    and that $\mathcal{G}(R) \cap H(K)$ is dense in~$H$.
    Then $\mathcal{G}$ satisfies (OPS) for codimension at least $\dim {\mathcal{G}_K}+1$.
\end{proposition}

We stress that $\mathcal{G}$ satisfying (OPS) for codimension at least~$\dim G_K + 1$ is not an empty statement, as the codimension of the closed subscheme $\mathcal{Z} \subseteq \mathcal{G}$ is taken in the generic fiber.
In fact, the codimension of $\mathcal{Z}_K$ in $\mathcal{G}_K$ being at least $\dim G_K + 1$ is equivalent to $\mathcal{Z}_K$ being empty, which is equivalent to
the image of $\mathcal{Z}$ in $\Spec R$ being a finite set of closed points by \cite[Corollaire~9.7.9]{EGAIV3}.
If $\mathcal{X}$ is a finite type $R$-scheme and~$\mathcal{Z}\subseteq \mathcal{X}$ is a closed subscheme whose image in $\Spec R$ is a finite set of closed points, we say that $\mathcal{Z}$ is a \emph{vertical} closed subscheme of $\mathcal{X}$.

To prove Proposition~\ref{Prop:HT_general}, we have to show that, given a vertical closed sub\-scheme~$\mathcal{Z}\subseteq \mathcal{G}$,
the non-emptiness of $(\mathcal{G}\setminus \mathcal{Z})(R)$ implies its non-strongly-thinness.
We start by showing the existence of a dense subgroup of integral points, to which we may then apply a fibration theorem for the weak Hilbert property of algebraic groups given in \cite{LugerProducts}.
For the existence of such a subgroup, we use an approximation argument
given by Hassett--Tschinkel \cite[Proposition~4.1]{HT01}, which they attribute to McKinnon:

\begin{lemma}\label{Lem:HT}
    Let $\mathcal{G}$ be a finite type group scheme over $R$ and
    let $\mathcal{Z}$ be a vertical closed subscheme of~$\mathcal{G}$.
    Then the following statements hold.
    \begin{enumerate}[label=(\arabic*)]
    	\item\label{Item:LemHT_coset}
    		For every point $p\in (\mathcal{G}\setminus \mathcal{Z})(R)$, there exists an integer $m\geq 1$ such that the inclusion $p\cdot \mathcal{G}(R)^m \subseteq (\mathcal{G}\setminus \mathcal{Z})(R)$ holds.
    	\item\label{Item:LemHT_dense}
    		Assume that $\mathcal{G}$ is connected and $\mathcal{G}(R)$ is dense.
    		If $(\mathcal{G}\setminus \mathcal{Z})(R)$ is non-empty, then it is dense.
    \end{enumerate} 
\end{lemma}
\begin{proof}
Since $\mathcal{Z}$ is vertical, its image in $\Spec R$ is finite.
Let $\mathfrak{p}_1,\ldots,\mathfrak{p}_n$ be the primes lying under~$\mathcal{Z}$.
Let $m_i $ be the order of the finite group $\mathcal{G}(k(\mathfrak{p}_i))$ and let $m$ denote the least common multiple of the $m_i$.
Note that, for every $g$ in $\mathcal{G}(R)^m\subseteq \mathcal{G}(R)$, the element~$p\cdot g$ lies in $(\mathcal{G} \setminus \mathcal{Z})(R)$. Indeed, for every $i=1,\ldots,n$, the image of $p\cdot g$ in~$\mathcal{G}(k(\mathfrak{p}_i))$ equals the image of $p$, as $g$ reduces to the identity element.
This concludes the proof of \ref{Item:LemHT_coset}.

To prove \ref{Item:LemHT_dense}, let $p\in (\mathcal{G}\setminus \mathcal{Z})(R)$. By \ref{Item:LemHT_coset} there exists an integer $m\geq 1$ such that $p\cdot \mathcal{G}(R)^m\subseteq (\mathcal{G}\setminus \mathcal{Z})(R)$.
Let $[m]\colon \mathcal{G}_K \to \mathcal{G}_K$ denote the $m$-power map.
By \cite[Lemma~4.4.12]{SpringerLinAlgGrps}, the tangent map induced by $[m]$ is the multiplication by $m$ and thus surjective.
Since $\mathcal{G}_K$ is smooth over $K$ \cite[II.§5, Corollaire 2.3]{DemazureGabriel}, it follows from
\cite[Lemma~3.3]{BartschJavanpeykar} that $[m]$ is dominant.
Thus, by the density of $\mathcal{G}(R)$, its image $\mathcal{G}(R)^m$ under $[m]$ is dense,
so $(\mathcal{G}\setminus \mathcal{Z})(R)$ contains the dense subset $p\cdot \mathcal{G}(R)^m$ and is thus dense.
\end{proof}

The proof of Proposition~\ref{Prop:HT_general} is now an application of the fibration results from \cite{LugerProducts}.

\begin{proof}[Proof of Proposition~\ref{Prop:HT_general}]
    Let $\mathcal{Z} \subseteq \mathcal{G}$ be a closed subscheme with $(\mathcal{G}\setminus \mathcal{Z})(R) \neq \emptyset$
    and such that $\mathcal{Z}_K$ is of codimension at least~$\dim {\mathcal{G}_K}+1$ in $\mathcal{G}_K$, i.e. $\mathcal{Z}_K = \emptyset$, so $\mathcal{Z}$ is a vertical closed subscheme of $\mathcal{G}$.
    Let $p\in (\mathcal{G}\setminus \mathcal{Z})(R)$. By part \ref{Item:LemHT_coset} of Lemma~\ref{Lem:HT}, we may choose $m\geq 1$ such that $p\cdot \mathcal{G}(R)^m\subseteq (\mathcal{G} \setminus \mathcal{Z})(R)$.
    To prove the claim, it suffices to show that $\Omega \coloneqq \mathcal{G}(R)^m$
    is not strongly thin is $\mathcal{G}_K$, since we may compose ramified covers of~$\mathcal{G}_K$ with translation by $p$.
    Note that $\Omega$ is dense in $\mathcal{G}_K$ by part \ref{Item:LemHT_dense} of Lemma~\ref{Lem:HT} and, similarly, $\Omega \cap H(K) \supseteq (\mathcal{G}(R)\cap H(K))^m$ is dense in $H$.
    Moreover, $\Omega$ is a subgroup of $\mathcal{G}(K)$ by the commutativity of $\mathcal{G}$.
    Thus, by \cite[Proposition~5.10]{LugerProducts}, there exists a finitely generated subgroup $\Omega'$ of $\Omega$ such that $\Omega'$ is dense in $\mathcal{G}$ and $\Omega' \cap H(K)$ is dense in $H$.
    By \cite[Proposition~5.9]{LugerProducts}, this implies that $\Omega'$, is not strongly thin in $\mathcal{G}_K$, and thus $\Omega$ is not strongly thin in $\mathcal{G}_K$.
\end{proof}

\begin{corollary}\label{cor_ht_new}
    Let $\mathcal{G}$ be either a linear connected commutative finite type group scheme over~$R$
    or an abelian scheme over $R$.
    If $\mathcal{G}_K$ is connected and $\mathcal{G}(R)$ is dense, then $\mathcal{G}$ satisfies (OPS) for codimension at least~$\dim {\mathcal{G}_K}+1$.
\end{corollary}
\begin{proof}
    Apply Proposition~\ref{Prop:HT_general} to the obvious exact sequence, in which $H=0$ if $\mathcal{G}$ is an abelian scheme, or $A=0$ if $\mathcal{G}$ is linear.
\end{proof}

\begin{corollary}\label{Cor:CommGroupsOPSE}
    Let $\mathcal{G}$ be a connected commutative finite type group scheme over $R$ with $\mathcal{G}_K$ connected and $\mathcal{G}(K)$ dense.
    Then $\mathcal{G}$ satisfies (OPSE) for codimension at least~$\dim {\mathcal{G}_K}+1$.
\end{corollary}
\begin{proof}
    Note that $\mathcal{G}_K$ is smooth over $K$ by \cite[II.§5, Corollaire 2.3]{DemazureGabriel}. Thus, by \cite[Theorem~1.1]{ConradCT}, there exists an exact sequence
    \[
        0 \longrightarrow H \longrightarrow \mathcal{G}_K \longrightarrow A \longrightarrow 0
    \]
    of connected algebraic groups over $K$ with $H$ linear and $A$ an abelian variety.
    By assumption, $\Omega \coloneqq \mathcal{G}(K)$ is dense in $\mathcal{G}_K$,
    and by \cite[Corollary~18.3]{BorelLinAlgGrp}, $H(K) = \Omega \cap H(K)$ is dense in $H$.
    Thus, by \cite[Proposition~5.10]{LugerProducts}, there exists a finitely generated dense subgroup $\Omega' \subseteq \Omega$ such that $\Omega' \cap H(K)$ is dense in $H$.
    Let $R' = \mathcal{O}_{K,S'}$ for $S'$ a finite set of finite places of $K$ containing $S$ such that the generators of $\Omega'$
    lie in $\mathcal{G}(R')$, which implies that $\Omega' \subseteq \mathcal{G}(R')$.
    Let $S''$ be any finite set of finite places of $K$ containing $S'$ and define $R'' = \mathcal{O}_{K,S''}$. Since $\mathcal{G}(R'') \supseteq \mathcal{G}(R')$, we see that
    $\mathcal{G}(R'') \supseteq \mathcal{G}(R')$ is dense in $\mathcal{G}_K$ and $\mathcal{G}(R'')\cap H(K)$ is dense in $H$, so $\mathcal{G}_{R''}$ satisfies (OPS) for codimension at least~$\dim \mathcal{G}_K + 1$ by Proposition~\ref{Prop:HT_general}.
\end{proof}

\begin{example}
    Let $\mathcal{G} = \Gmn{R}^4$ and let $\mathcal{Z}\subseteq \mathcal{G}$ be a closed subscheme.
    \begin{enumerate}[label=(\arabic*)]
    \item By Corollary~\ref{Cor:CommGroupsOPSE}, each factor $\Gmn{R}$ satisfies (OPSE) for codimension at least~$2$.
    Thus, if $\mathcal{Z}_K \subseteq \mathcal{G}_K$ is of codimension at least~$2$, then by Theorem~\ref{Thm:Products_OPSE}, there exists a finite set of finite places $S'$ containing $S$ such that
    $(\mathcal{G}\setminus \mathcal{Z})(\mathcal{O}_{K,S'})$ is not strongly thin in $\mathcal{G}_K$.
    \item Assume that $(\mathcal{G}\setminus \mathcal{Z})(R)$ is dense and that $\mathcal{Z}_K \subseteq \mathcal{G}_K$ is of codimension at least~$3$.
    Write $\mathcal{G} \cong \Gmn{R}^2 \times_R \Gmn{R}^2$.
    Then $\Gmn{R}^2(R)$ is dense in $\Gmn{K}^2$, so both factors satisfy (OPS) for codimension at least~$3$ by Corollary~\ref{cor_ht_new}. Thus, $(\mathcal{G}\setminus \mathcal{Z})(R)$ is already not strongly thin in $\mathcal{G}_K$ by Theorem~\ref{Thm:Products_OPS}.
    \end{enumerate}
\end{example}

\section{Proof of Theorem~\ref{Thm:WHP_Punctured_LinAlgGrps}}

In this section, we prove Theorem~\ref{Thm:WHP_Punctured_LinAlgGrps}, first for semisimple quasi-split linear algebraic groups and unipotent groups, and then in general by reducing to these cases.
We first collect some results pertaining to strong approximation and the usual Hilbert property.

\subsection{Preliminaries on strong approximation and the Hilbert property}

Recall that an integral variety over a field $k$ satisfies the Hilbert property if its set of~$k$-rational points is not thin, and that a strongly thin set in a normal $k$-variety is thin:

\begin{definition} 
Let $X$ be an integral variety over a field  $k$.
A subset $\Omega\subseteq X(k)$  is \emph{thin} in $X$ if there is an integer $n\geq 0$ and a finite collection $(\pi_i \colon Y_i\to X)_{i=1}^n$ of finite covers of degree at least~$2$ such that  
\[
    \Omega \setminus \bigcup_{i=1}^n \pi_i (Y_i(k))
\]
is not dense in $X$.
\end{definition}

If $U\subseteq X$ is a dense open and $\Omega \subseteq U(k)$,
then $\Omega$ is thin in $U$ if and only if $\Omega$ is thin in $X$.
Note that the ``if'' part of the same statement might fail if we replace ``thin'' by ``strongly thin''.

A classical argument due to Ekedahl (see \cite[Theorem~3.5.7]{SerreTopicsGalois}) shows that the Hilbert property is implied by the property ``weak weak approximation'' (see \cite[Definition~3.5.6]{SerreTopicsGalois}).
In this subsection we illustrate that, similarly, the integral Hilbert property is implied by ``strong approximation'',
which is quickly derived from work due to Moerman \cite{Moerman} and Nakahara--Streeter \cite{NakaharaStreeter},
and relate this to the properties~(OPS) and (OPSE).

Recall that the ring of \emph{adeles} of $K$ is the subring $\adeles_K \subseteq \prod_v K_v$,
where the product is taken over all finite places $v$ of $K$,
whose elements $(\alpha_v)$ satisfy $\alpha_v \in \mathcal{O}_{K,v}$ for all except finitely many $v$.
For a finite set of finite places $T$ of $K$,
the ring $\adeles_K^T$ of \emph{adeles away from $T$} is the subring
$\adeles_K^T = \prod_{v\notin T} \mathcal{O}_{K,v} \times \prod_{v\in T} K_v$ of $\adeles_K$.
Moreover, one defines analogously the ring $\adeles_{K,S}$ of \emph{$S$-adeles} and the ring $\adeles_{K,S}^T$ of \emph{$S$-adeles away from~$T$} as subrings of $\prod_{v\notin S} K_v$.
Note that $K$ embeds diagonally into $\adeles_K$ and that~$\adeles_K$, resp.~$\adeles_{K,S}$, projects onto $\adeles_K^T$, resp. $\adeles_{K,S}^T$.
In particular, given a variety $X$ over $K$, there are natural maps
$X(K) \to X(\adeles_K) \to X(\adeles_K^T)$ to the adelic points of $X$ and the adelic points of $X$ away from~$T$. These spaces carry the adelic topology;
see \cite{ConradAdeles} for more details on this and for the central equations of the following paragraph.

A variety $X$ over $K$ is said to satisfy \emph{strong approximation off $T$} if
the diagonal image of
$X(K)$ is dense in $X(\adeles_K^T)$, see \cite[Definition~12.2.7]{CTS_Brauer}.
If $\mathcal{X}$ is a model of $X$ over $\mathcal{O}_{K,S}$,
it then follows that
$\mathcal{X}(\mathcal{O}_{K,S\cup T})$ is dense in
\[
    \mathcal{X}(\adeles_{K,S\cup T}^T)
    = \mathcal{X}(\adeles_{K,S}^T)
    = \prod_{v \notin S\cup T} \mathcal{X}(\mathcal{O}_{K,v}) \times \prod_{v\in S\setminus T}\mathcal{X}(K_v)
    \subseteq X(\adeles_K^T)
.\]
For details, see the proof of Proposition~\ref{Prop:SA_implies_HP}.
If $X$ is smooth,
this property of $\mathcal{X}$ is referred to as \emph{integral strong approximation off $T$} 
in \cite[Definition~2.21]{MitNakStr}; we will use the term integral strong approximation without the smoothness assumption.
If~$S$ contains $T$ and $\mathcal{X}$ is normal, integral strong approximation off $T$ is a special case of \emph{Campana weak weak approximation}, as defined in \cite[Definition~2.7]{NakaharaStreeter}.
Note that the definition given by Nakahara--Streeter requires $\mathcal{X}$ to be flat over $\mathcal{O}_{K,S}$, which is satisfied here by \cite[Proposition~4.3.9]{LiuAGAC} since $\mathcal{X}$ is dominant over $\mathcal{O}_{K,S}$ and $\mathcal{O}_{K,S}$ is a Dedekind domain.

If $\mathcal{X}$ is normal, it is shown in \cite[Theorem~1.1]{NakaharaStreeter} that Campana weak weak approximation implies the \emph{Campana Hilbert property}, which is the usual Hilbert property in our case.
This theorem of Nakahara--Streeter has been further generalized by Moerman \cite[Theorem~1.1]{Moerman}, where normality of $\mathcal{X}$ is not required, which gives us the following proposition.

\begin{lemma}\label{Lem:ISA_implies_HP}
    Let $T$ be a finite set of finite places of $K$, let $\mathcal{X}$ be an arithmetic scheme over $\mathcal{O}_{K,S}$, and assume that $\mathcal{X}$ satisfies integral strong approximation off $T$, i.e.,
    $\mathcal{X}(\mathcal{O}_{K,S\cup T})$ is dense in $\mathcal{X}(\adeles_{K,S \cup T}^T)$.
    If $\mathcal{X}(\mathcal{O}_{K,S\cup T})$ is not empty, then it is not thin in $\mathcal{X}_K$.
\end{lemma}
\begin{proof}
    To see that this is a special case of \cite[Theorem~1.1]{Moerman},
    we write the integral points of $\mathcal{X}$ using the notation of $\mathcal{M}$-points introduced in \cite[Definition~3.12]{Moerman}.
    By the universal property of fiber products, we may replace $\mathcal{X}$ by $\mathcal{X}_{\mathcal{O}_{K,S\cup T}}$ to assume that $S$ contains $T$.
    
    Let $\mathcal{X} \to \overline{\mathcal{X}}$ be an open immersion of arithmetic $\mathcal{O}_{K,S}$-schemes with $\overline{\mathcal{X}}$ proper over $\mathcal{O}_{K,S}$
    and define $\mathcal{D} = \overline{\mathcal{X}} \setminus \mathcal{X}$.
    Let $X, \overline{X}, D$ denote the generic fiber of
    $\mathcal{X}, \overline{\mathcal{X}}, \mathcal{D}$, respectively.
    We consider the set of multiplicities $\mathfrak{M} = \{(0)\}$ and define
    $\mathcal{M} = ((\mathcal{D}), \mathfrak{M})$ and
    $M = ((D), \mathfrak{M})$.
    Then $\overline{X}$ is a proper variety over the PF field $(K, \Spec \mathcal{O}_K)$
    and the pair $(\overline{\mathcal{X}}, \mathcal{M})$ is an integral model of the pair
    $(\overline{X}, M)$ over $\Spec \mathcal{O}_{K,S}$ (see~\cite[Definition~3.1 and Definition~3.4]{Moerman}).
    Moreover, for the pair $(\overline{\mathcal{X}}, \mathcal{M})$, the set of
    $\mathcal{M}$-integral points is $\mathcal{X}(\mathcal{O}_{K,S})$ \cite[\S 3.3, Example~(1)]{Moerman}
    and the set of integral adelic~$\mathcal{M}$-points is $\mathcal{X}(\adeles_{K,S \cup T}^T)$ \cite[Example~4.4]{Moerman}.

    Our assumption now says that
    $(\overline{\mathcal{X}}, \mathcal{M})$ admits an $\mathcal{M}$-integral point and that
    its~$\mathcal{M}$-integral points are dense in its integral adelic $\mathcal{M}$-points, i.e.,
    $(\overline{\mathcal{X}}, \mathcal{M})$ satisfies integral~$\mathcal{M}$-approximation off $T$ \cite[Definition~4.8]{Moerman}.
    Therefore, by \cite[Theorem~1.1]{Moerman}, we obtain that $\mathcal{X}(\mathcal{O}_{K,S})$ is not thin in $\overline{X}$, and thus not in $X$, as required.
\end{proof}

\begin{proposition}\label{Prop:SA_implies_HP}
    Let $T$ be a finite set of finite places of $K$ contained in $S$,
    let $X$ be an integral $K$-variety satisfying strong approximation off $T$,
    and let $\mathcal{X}$ be a model of~$X$ over $\mathcal{O}_{K,S}$.
    If $\mathcal{X}(\mathcal{O}_{K,S})$ is non-empty, then it is not thin in $X$.
\end{proposition}
\begin{proof}
    As $\mathcal{X}(\mathcal{O}_{K,S})$ is non-empty,
    $\mathcal{X}(\adeles_{K,S}^T)$ is a non-empty open subset of $X(\adeles_K^T)$.
    By strong approximation off $T$, the diagonal image of $X(K)$ inside $X(\adeles_K^T)$ is dense.
    Therefore, the intersection $\mathcal{X}(\adeles_{K,S}^T) \cap X(K)$ is dense in $\mathcal{X}(\adeles_{K,S}^T)$.
    Since every $K$-point of $X$ inside $\mathcal{X}(\adeles_{K,S}^T)$ is $S$-integral, this shows that
    the diagonal image of $\mathcal{X}(\mathcal{O}_{K,S})$ in~$\mathcal{X}(\adeles_{K,S}^T)$ is dense.
    Thus, it follows from Lemma~\ref{Lem:ISA_implies_HP} that $\mathcal{X}(\mathcal{O}_{K,S})$ is not thin in $X$, as claimed.
\end{proof}

If $X$ is a variety over $K$ and $X$ satisfies strong approximation off a finite set of finite places $T$ of $K$, it was asked by Wittenberg \cite[Question~2.11]{WittenbergRationalPoints} if $X\setminus Z$ satisfies strong approximation off $T$ for every closed subscheme $Z\subseteq X$ of codimension at least~$2$.
If this condition is satisfied for a given normal integral variety $X$ and $\mathcal{X}$ is a model of $X$ with an integral point over $\mathcal{O}_{K,S}$,
then Proposition~\ref{Prop:SA_implies_HP} implies that $\mathcal{X}$ satisfies (OPS) for codimension at least~$2$.

For $T \neq \emptyset$, a positive answer to Wittenberg's question was given
for $X=\mathbb{A}_K^n$ by Cao--Xu \cite[Proposition~3.6]{CaoXu} and Wei \cite[Lemma~2.1]{Wei},
and for $X$ a connected semi-simple simply connected quasi-split linear algebraic group by Cao--Liang--Xu~\cite[Theorem~1.1]{CaoLiangXu}.
The latter result has been further extended by Cao--Huang~\cite[Theorem~1.4]{CaoHuang} to semi-simple simply connected almost simple algebraic groups which are isotropic.
As a direct corollary of these results on strong approximation and the work of Moerman and Nakahara--Streeter, we obtain the following.

\begin{corollary}\label{Cor:CLX_models_OPS}
    Let $\mathcal{G}$ be an arithmetic scheme over $\mathcal{O}_{K,S}$ such that $\mathcal{G}_K$
    is a quasi-split semi-simple simply connected linear algebraic group over $K$ or a unipotent group over $K$.
    If $S\neq \emptyset$ and $\mathcal{G}(\mathcal{O}_{K,S}) \neq \emptyset$, then $\mathcal{G}$ satisfies (OPS) for codimension at least~$2$.
\end{corollary}
\begin{proof}
    Let $\mathcal{Z}\subseteq \mathcal{G}$ be a closed subscheme such that $\mathcal{Z}_K \subseteq \mathcal{G}_K$ is of codimension at least~$2$.
    Define $\mathcal{X} = \mathcal{G}\setminus \mathcal{Z}$ and $X=\mathcal{X}_K$ and assume that $\mathcal{X}(\mathcal{O}_{K,S}) \neq \emptyset$.
    If $\mathcal{G}_K$ is a semi-simple simply connected quasi-split linear algebraic group over $K$, then $X$ satisfies strong approximation off $S$ by \cite[Theorem~1.1]{CaoLiangXu}.
    If $\mathcal{G}_K$ is a unipotent group over $K$, then it is isomorphic (as a variety) to $\mathbb{A}_K^n$, and thus $X$ satisfies strong approximation off $S$ by \cite[Proposition~3.6]{CaoXu} (see also \cite[Lemma~2.1]{Wei}).
    Therefore, it follows from Proposition~\ref{Prop:SA_implies_HP} (taking $T=S$) that $\mathcal{X}(\mathcal{O}_{K,S})$ is not thin in $X$ and thus not thin in~$\mathcal{G}_K$.
\end{proof}

\subsection{Application to linear algebraic groups}

Using universal coverings and a Chevalley--Weil type lifting arguments for integral points, we can drop the simply-connectedness assumption in Corollary~\ref{Cor:CLX_models_OPS} by replacing the property (OPS) with (OPSE).

\begin{proposition}\label{Prop:QSplit_SSimple_Unipotent_OPSE}
    Let $\mathcal{G}$ be an arithmetic scheme over $\mathcal{O}_{K,S}$ such that $\mathcal{G}_K$ is a semi-simple quasi-split linear algebraic group over $K$ or a unipotent group over $K$.
    Then $\mathcal{G}$ satisfies (OPSE) for codimension at least~$2$.
\end{proposition}
\begin{proof}
    Let $G = \mathcal{G}_K$. Assume first that $G$ is unipotent.
    Let $S' \neq \emptyset$ be a finite set of finite places of $K$ containing $S$ such that $\mathcal{G}(\mathcal{O}_{K,S'})\neq \emptyset$.
    Then, by Corollary~\ref{Cor:CLX_models_OPS}, the base change $\mathcal{G}_{\mathcal{O}_{K,S''}}$ satisfies (OPS) for codimension at least~$2$ for every finite set of finite places $S''$ of $K$ containing $S'$, i.e., $\mathcal{G}$ satisfies (OPSE) for codimension at least~$2$.

    Now assume that $G$ is semi-simple quasi-split.
    By \cite[Corollary~A.4.11]{CGP_PseudoReductive}, there exists an isogeny $\lambda \colon G' \to G$ with $G'$ a connected semi-simple simply connected linear algebraic group over~$K$.
    Enlarging $S$ if necessary, there exists a model $\mathcal{G}'$ of~$G'$ over~$R = \mathcal{O}_{K,S}$ such that $\lambda$ extends to a finite étale morphism $\bar \lambda\colon \mathcal{G}' \to \mathcal{G}$.
    Let $S'$ be any finite set of finite places of $K$ containing $S$.
    To prove the claim, we may replace~$S$ by $S'$, and $\mathcal{G}$ and $\mathcal{G}'$ by $\mathcal{G}_{\mathcal{O}_{K,S'}}$ and $\mathcal{G}_{\mathcal{O}_{K,S'}}$ accordingly, and show that $\mathcal{G}$ satisfies (OPS) for codimension at least~$2$.
    
    By \cite[Theorem~3.2]{LugerProducts}, there exists a
    finite étale morphism $\Spec T \to \Spec R$ such that $\mathcal{G}(R)$ is contained in $\bar \lambda(\mathcal{G}'(T))$.
    Let $L$ denote the fraction field of $T$.
    We modify a standard application of Weil restriction of scalars, as used for example in the proofs of~\cite[Proposition~1.5]{CZHP} or \cite[Theorem~3.9]{CDJLZ}, to integral points to replace the $T$-integral points of $\mathcal{G}'$ by $R$-integral points of finitely many twists of $\mathcal{G}'$.
    Consider the induced morphism of Weil restrictions
    \[
        \Res_{T/R} (\bar \lambda_T) \colon
        \Res_{T/R} \mathcal{G}'_T \to \Res_{T/R} \mathcal{G}_T
    .\]
    Since Weil restrictions are compatible with fiber products by the discussion following Lemma~1 of \cite[\S7.6]{BoschNeronModels}, we have that
    \[
        (\Res_{T/R} (\bar \lambda_T))_K = 
        \Res_{L/K} (\lambda_L) \colon
        \Res_{L/K} G'_L \to \Res_{L/K} G_L
    .\]
    Let $\Delta\colon \mathcal{G} \to \Res_{T/R} \mathcal{G}_T$ denote the diagonal embedding.
    Let $\mathcal{G}_1, \ldots, \mathcal{G}_d$ be the connected components of $\mathcal{G} \times_{\Delta, \Res_{T/R} \mathcal{G}_T} \Res_{T/R} \mathcal{G}'_T$ with $\mathcal{G}_i(R)\neq \emptyset$ and define~$G_i = \mathcal{G}_{i,K}$.
    Note that the natural morphisms $\lambda_i \colon G_i \to G$ are finite étale
    and extend to the natural morphisms $\bar\lambda_i \colon \mathcal{G}_i \to \mathcal{G}$,
    so in particular, we have $\lambda_i(\mathcal{G}_i(R)) \subseteq \mathcal{G}(R)$.
    Since $G_i$ is normal connected with a $K$-rational point, it is geometrically integral.
    Moreover, over the algebraic closure $\overline{K}$ of $K$, each $G_i$ is isomorphic to $G'_{\overline{K}}$, with the base change of~$\lambda_i$ over $\overline{K}$ being identified with $\lambda_{\overline{K}}$. In particular, the $G_i$ are connected semi-simple simply connected linear algebraic groups over $K$.
    Quasi-splitness is invariant under isogeny, since the image of a Borel subgroup under a surjective homomorphism of linear algebraic groups is a Borel subgroup, see \cite[Corollary~6.2.8]{SpringerLinAlgGrps}.
    In particular,
    since $G$ is quasi-split, so are the $G_i$.
    Thus, by Corollary~\ref{Cor:CLX_models_OPS}, all $\mathcal{G}_i$ satisfy (OPS) for codimension at least~$2$.

    Let $\mathcal{Z}\subseteq \mathcal{G}$ be a closed subscheme with $\mathcal{Z}_K \subseteq G$ of codimension at least~$2$ and such that $(\mathcal{G}\setminus\mathcal{Z})(R) \neq \emptyset$;
    we have to prove that $(\mathcal{G}\setminus\mathcal{Z})(R)$ is not strongly thin in~$G$.
    Let~$p\in (\mathcal{G}\setminus\mathcal{Z})(R)$.
    Since $\mathcal{G}(R)$ is contained in $\lambda(\mathcal{G}'(T))$, by the defining properties of Weil restriction, $p$ lifts to an $R$-integral point of $\Res_{T/R} \mathcal{G}'_T$, and thus to a point~$p' \in \mathcal{G}_i(R)$ for some $i$.
    Let $\mathcal{Z}'$ denote the inverse image of $\mathcal{Z}$ in $\mathcal{G}_i$.
    By \cite[Corollaire~6.1.4]{EGAIV2}, $\mathcal{Z}'_K \subseteq G_{i}$ is of codimension at least~$2$. Since $p' \in (\mathcal{G}_i \setminus \mathcal{Z}')(R)$ and $\mathcal{G}_i$ satisfies (OPS) for codimension at least~$2$, the set $(\mathcal{G}_i \setminus \mathcal{Z}')(R)$ is not strongly thin in $G_i$.
    Moreover, its image under the finite étale morphism $\lambda_i$ is contained in $(\mathcal{G}\setminus\mathcal{Z})(R)$.
    By \cite[Lemma~5.3(b)]{BSFP23BaseChange}, this implies that $(\mathcal{G}\setminus\mathcal{Z})(R)$ is not strongly thin in $G$, as required.
\end{proof}

For the sake of clarity, we point out that Proposition~\ref{Prop:QSplit_SSimple_Unipotent_OPSE} says that connected semi-simple quasi-split linear algebraic groups and unipotent groups over $K$ satisfy the ``punctured integral Hilbert property'', i.e., part \ref{Item:Conj_XminusZ_IWHP} of Conjecture~\ref{conj}, even without enlarging $K$.
Combining this with Theorem~\ref{Thm:Products_OPSE} and Corollary~\ref{Cor:CommGroupsOPSE}, we obtain the same result for products of those groups and $1$-dimensional connected commutative algebraic groups with a dense set of rational points.

\begin{corollary}\label{Cor:Prod_1-DimGroups}
    Let $G = \prod_{i=1}^N G_i$ be an algebraic group over $K$ with each $G_i$ either isomorphic to
    $\mathbb{G}_{\mathrm{a}}$, $\Gm$, a connected semi-simple quasi-split linear algebraic group, a unipotent group, or an elliptic curve over $K$.
    Let $Z \subseteq G$ be a closed subscheme of codimension at least~$2$.
    If $G(K)$ is dense, then there exists a finite set $T$ of finite places of $K$ such that $G \setminus Z$ satisfies the weak Hilbert property over $\mathcal{O}_{K,T}$.
    In particular, $G$ satisfies part \ref{Item:Conj_XminusZ_IWHP}, and thus also \ref{Item:Conj_X_denseIntegralPoints}-\ref{Item:Conj_X_IWHP}, of Conjecture~\ref{conj}.
\end{corollary}
\begin{proof}
    Let $T$ be a finite set of finite places of $K$ such that, for each $i$, there exists
    a normal quasi-projective model $\mathcal{G}_i$ of $G_i$ over $\mathcal{O}_{K,T}$,
    and $\mathcal{G}_i$ is a commutative group scheme over $\mathcal{O}_{K,T}$ for all~$i$ for which $\dim G_i = 1$.
    For every index $i$ for which $G_i$ is not of dimension $1$, it is semi-simple quasi-split or unipotent.
    Let $R' = \mathcal{O}_{K,T}$,
    consider the model $\mathcal{G} = \mathcal{G}_i \times_{R'} \ldots \times_{R'} \mathcal{G}_N$ of $G$ over $R'$,
    and let $\mathcal{Z} \subseteq \mathcal{G}$ be the closure of $Z$ in~$\mathcal{G}$.
    By Corollary~\ref{Cor:CommGroupsOPSE} and Proposition~\ref{Prop:QSplit_SSimple_Unipotent_OPSE},
    each $\mathcal{G}_i$ satisfies (OPSE) for codimension at least~$2$,
    so the claim follows from Theorem~\ref{Thm:Products_OPSE}.
\end{proof}

We use the following descent argument to reduce the proof of Theorem~\ref{Thm:WHP_Punctured_LinAlgGrps}
to products of algebraic groups that satisfy (OPSE).

\begin{lemma}\label{Lem:PPWHP_descends_smooth_proper}
    Let $\varphi\colon X' \to X$ be a smooth proper morphism of normal integral varieties over $K$ and let $c \geq 1$ be an integer.
    Assume that, for every closed subscheme~$Z' \subseteq X'$ of codimension at least~$c$, there is a number field $L$ and a finite set of finite places $T$ of $L$ such that $(X'\setminus Z')_L$ satisfies the  weak Hilbert property over~$\mathcal{O}_{L,T}$.
    Then, for every closed subscheme $Z \subseteq X$ of codimension at least~$c$, there is a number field $L$ and a finite set of finite places $T$ of $L$ such that $(X\setminus Z)_L$ satisfies the  weak Hilbert property over $\mathcal{O}_{L,T}$.
\end{lemma}
\begin{proof}
    Let $Z \subseteq X$ be a closed subscheme of codimension at least~$c$.
    Let $Z' = \varphi^{-1}(Z)$. By \cite[Corollaire~6.1.4]{EGAIV2},
    $Z' \subseteq X'$ has codimension at least~$c$.
    By the assumption on~$X'$, there is a number field $L$, a finite set of finite places $T$ of $L$,
    and a model $\mathcal{U}'$ of $(X'\setminus Z')_L$ over $\mathcal{O}_{L,T}$ such that
    $\mathcal{U}'(\mathcal{O}_{L,T})$ is not strongly thin in $\mathcal{U}'_L$.
    Enlarging $T$ if necessary, there exists a model $\mathcal{U}$ of $(X\setminus Z)_L$ over $\mathcal{O}_{L,T}$
    such that the smooth proper morphism $\varphi' \colon \mathcal{U}_L' \to \mathcal{U}_L$ induced by $\varphi$ extends to a morphism $\mathcal{U}' \to \mathcal{U}$ over $\mathcal{O}_{L,T}$.
    In particular, the image $\varphi'(\mathcal{U}'(\mathcal{O}_{L,T}))$ is contained in
    $\mathcal{U}(\mathcal{O}_{L,T})$.
    Since $\mathcal{U}'(\mathcal{O}_{L,T})$ is not strongly thin in $\mathcal{U}'_L$, it follows from \cite[Lemma~5.3(b)]{BSFP23BaseChange} that the image $\varphi'(\mathcal{U}'(\mathcal{O}_{L,T}))$,
    and thus $\mathcal{U}(\mathcal{O}_{L,T})$, is not strongly thin in $\mathcal{U}_L = (X \setminus Z)_L$.
    This shows that $(X\setminus Z)_L$ satisfies the  weak Hilbert property over $\mathcal{O}_{L,T}$.
\end{proof}

Note that Lemma~\ref{Lem:PPWHP_descends_smooth_proper} says in particular that condition~\ref{Item:Conj_XminusZ_IWHP} of Conjecture~\ref{conj}
descends along smooth proper morphisms.
The same also holds for smooth surjective morphisms with geometrically integral generic fiber by \cite[Lemma~5.3(a)]{BSFP23BaseChange}.

\begin{proof}[Proof of Theorem~\ref{Thm:WHP_Punctured_LinAlgGrps}]
By Mostow's Theorem \cite[\S VIII, Theorem~4.3]{HochschildAlgGrps}, there exists a unipotent group $U$ and a reductive group $P$ over $K$ such that $G$ is the semi-direct product of $U$ and $P$.
In particular, $G$ is isomorphic to $U \times P$ as a $K$-variety.

Let $\Rad(P)$ denote the radical of $P$ and let $P^{\der}$ denote the derived subgroup of $P$.
By \cite[Exposé~XXII, Proposition~6.2.4]{SGA3}, there exists an isogeny $\Rad(P) \times P^{\der} \to P$, which is in particular a finite étale cover.
Since $\Rad(P)$ is a torus, there exists a finite field extension $L/K$ and an integer $n\geq 0$ such that $\Rad(P)_L \cong \Gmn{L}^n$.
The derived subgroup $P^{\der}$ of $P$ is semi-simple,
so by \cite[Corollary~A.4.11]{CGP_PseudoReductive} there exists a connected semi-simple simply connected linear algebraic group $G'$ over $K$ and a finite étale cover~$G' \to P^{\der}$.
Enlarging $L$ if necessary, we may further assume that $G'_L$ is split.
Define $G'' \coloneqq U_L \times \Gmn{L}^n \times G'_L$
and note that we have a finite étale morphism~$G'' \to G_L$.

We claim that the theorem holds for $G''$. Let $Z'' \subseteq G''$ be a closed subscheme of codimension at least~$2$.
Let $T \neq \emptyset$ be a finite set of finite places of $L$ such that there exist normal quasi-projective models $\mathcal{U}$ of $U_L$ and $\mathcal{G'}$ of $G'$ over $R' \coloneqq \mathcal{O}_{L,T}$.
Let~$\mathcal{Z}''$ denote the closure of $Z''$ in the model
$\mathcal{G}'' \coloneqq \mathcal{U} \times_{R'} \Gmn{\mathcal{O}_{L,T}}^n \times_{R'} \mathcal{G'}$ of $G''$ over~$R'$.
Since~$\Gmn{R}$ satisfies (OPSE) for codimension at least~$2$ by Corollary~\ref{Cor:CommGroupsOPSE} and~$\mathcal{U}$ and $\mathcal{G}''$ satisfy (OPSE) for codimension at least~$2$ by Corollary~\ref{Cor:CLX_models_OPS}, it follows from Theorem~\ref{Thm:Products_OPSE} that
there exists a finite set of finite places $T'$ of $L$ containing $T$ such that $(\mathcal{G}''\setminus \mathcal{Z}'')(\mathcal{O}_{L,T'})$ is not strongly thin in $G''$,
i.e., the theorem holds for $G''$.
By Lemma~\ref{Lem:PPWHP_descends_smooth_proper}, this implies that it also holds for $G$.
\end{proof}

\bibliography{references}{}
\bibliographystyle{alpha}

\end{document}